\documentclass[titlepage,11pt]{article}
\usepackage{latexsym}
\usepackage{amsfonts}
\usepackage{amsmath}
\usepackage{tikz}
\usepackage{textcomp}
\usetikzlibrary{graphs}
\usepackage{float}

\newcommand\blackslug{\hbox{\hskip 1pt \vrule width 4pt height 8pt depth 1.5pt
        \hskip 1pt}}
\newcommand\bbox{\hfill \quad \blackslug \bigbreak}
\def\dd{\hbox{-}}
\def\cc{\hbox{-}\cdots\hbox{-}}
\def\ll{,\ldots,}

\title{Excluding the fork and antifork\footnote{Also available as: Maria Chudnovsky, Linda Cook, Paul Seymour, Excluding the fork and antifork, Discrete Mathematics, Volume 343, Issue 5, 2020.}}
\author{Maria Chudnovsky\thanks{This material is based upon work supported in part by the U. S. Army
Research Office under grant   number W911NF-16-1-0404, and by  NSF grant DMS-1763817.},
Linda Cook\thanks{ Partially supported by NSF Award 1514606.}, and 
Paul Seymour\thanks{Partially supported by AFOSR grant
A9550-19-1-0187 and NSF grant  DMS-1800053.}\\
Princeton University, Princeton, NJ 08544}

\date{February 22, 2019; revised November 17, 2019}

\newtheorem{thm}{}[section]

\newcommand{\Proof}{\noindent{\bf Proof.}\ \ }

\begin{document}
\maketitle
\begin{abstract}
The {\em fork} is the tree obtained from the claw $K_{1,3}$ by subdividing one of its edges once, and
the {\em antifork} is its complement graph. We give a complete description of all graphs that do not
contain the fork or antifork as induced subgraphs.
\end{abstract}

\section{Introduction}
\begin{figure}[H]
\centering

\begin{tikzpicture}[scale=0.8,auto=left]
\tikzstyle{every node}=[inner sep=1.5pt, fill=black,circle,draw]

\node (v0) at (0,0) {};
\node (v1) at (0,1) {};
\node (v2) at (0,2) {};
\node (v3) at (-.6,3) {};
\node (v4) at (.6,3) {};

\draw (v0) -- (v1);
\draw (v1) -- (v2);
\draw (v2) -- (v3);
\draw (v2) -- (v4);

\node (u0) at (6,0) {};
\node (u1) at (6,1) {};
\node (u2) at (6-.6,2) {};
\node (u3) at (6+.6,2) {};
\node (u4) at (6,3) {};

\draw (u0) -- (u1);
\draw (u1) -- (u2);
\draw (u1) -- (u3);
\draw (u2) -- (u3);
\draw (u2) -- (u4);
\draw (u3) -- (u4);

\end{tikzpicture}

\caption{The fork and the antifork.} \label{fig:fork}
\end{figure}
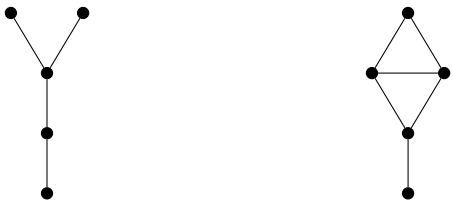

Graphs in this paper are finite, and without loops or parallel edges. The {\em fork} is the tree obtained from
a four-vertex path by adding a vertex adjacent to the second vertex of the path. The {\em antifork}
is its complement graph; thus, the antifork is obtained from a four-vertex path by adding one more vertex adjacent to the first three vertices
of the path. (See figure \ref{fig:fork}.)
Let us say $G$ is {\em uncluttered} if no induced subgraph of $G$ is a fork or antifork.
Our goal in this paper is to give a complete description of all 
uncluttered graphs. The line graph of a triangle-free graph is uncluttered, and so is its complement; and our main theorem says
that every uncluttered graph can be obtained by piecing together line graphs of triangle-free graphs and their complements.
Before we can make a precise statement we need to explain the ``piecing together'' process, and that is the content of the next section.
We state our main result in \ref{mainthm}.

This question was motivated by discussions with T. Karthick, who told us several results about colouring
graphs not containing forks and/or antiforks~\cite{alekseev, randerath,lozin} (note, however, that in \cite{randerath} 
``fork'' has a different meaning from its meaning here). In particular, he asked for the best possible 
``$\chi$-bounding function'' for uncluttered graphs. The answer follows from our main result \ref{mainthm},
and we give a proof in the final section.
\begin{thm}\label{chibounded}
For every uncluttered graph, its chromatic number is at most twice its clique number. 
\end{thm}
This is asymptotically best possible,
since if $H$ is a triangle-free graph with largest stable set of cardinality $k$, then the
complement of the line graph of $H$ is uncluttered, with clique number at most $|V(H)|/2$ and with chromatic number
$|V(H)|-k$; and we can choose $H$ and $k$ with $(|V(H)|-k)/|V(H)|$ arbitrarily close to $1$.

\section{Some safe operations}

There are several ways to make larger uncluttered graphs from smaller ones. The most obvious is: if $G_1, G_2$
are both uncluttered, then so is their disjoint union. The {\em complete join} of $G_1,G_2$ is obtained from their disjoint union
by adding edges between every vertex of $G_1$ and every vertex of $G_2$. Since the complement of an uncluttered graph is also 
uncluttered, and the complete join of $G_1,G_2$ is the complement of the disjoint union of $\overline{G_1},\overline{G_2}$,
it follows that if $G_1, G_2$
are both uncluttered, then so is their complete join.
Let us say $G$ is {\em anticonnected} if its complement graph $\overline{G}$
is connected. 

If $X,Y\subseteq V(G)$ are disjoint, we say $X$ is {\em complete} to $Y$ if every vertex in $X$ is adjacent to 
every vertex in $Y$; and $X$ is {\em anticomplete} to $Y$ if there are no edges between $X$ and $Y$.
If $v\in V(G)$, we say ``$v$ is complete to $Y$'' meaning that $\{v\}$ is complete to $Y$, and so on.

If $X\subseteq V(G)$ we denote by $G[X]$
the subgraph induced on $X$.
Two (distinct) vertices $u,v$ of $G$ are {\em twins} if $u,v$ have the same neighbours in $V(G)\setminus \{u,v\}$;
($u,v$ may or may not be adjacent).
A {\em homogeneous set} in $G$ means a set $X\subseteq V(G)$ such that every vertex of $G$ not in $X$ is either
complete or anticomplete to $X$; and the homogeneous set is {\em nontrivial} if $|X|\ge 2$ and $X\ne V(G)$.
If $X$ is a homogeneous set in $G$, let $Y$ be the set of vertices in $V(G)\setminus X$ that are complete to $X$.
Take a new vertex $v$, and let $H$ be the graph formed from $G\setminus X$ by adding the vertex $v$ and making $v$
adjacent to the vertices in $Y$. Then we say $G$ is obtained from $H$
by {\em substituting $G[X]$ for $v$}.

Say a vertex $v$ is {\em simplicial} if the set of its neighbours is a clique. If $v$ is a simplicial vertex
in an uncluttered graph $G$, then we may substitute a complete graph for $v$, and the new graph we obtain is also
uncluttered. Consequently, if $G$ has adjacent simplicial twins, then
$G$ can be obtained from a smaller graph by the operation just described. 
A vertex $v$ is {\em antisimplicial} if it is
simplicial in the complement graph, that is, if the set of all vertices nonadjacent to $v$ is a stable set.

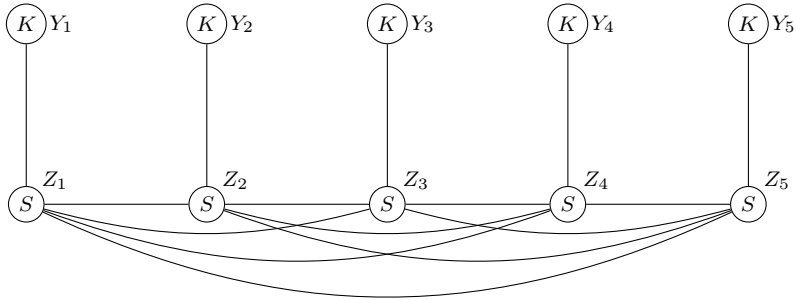
\begin{figure}[H]
\centering

\begin{tikzpicture}[scale=0.8,auto=left]
\tikzstyle{every node}=[inner sep=2pt, fill=white,circle,draw]

{\scriptsize
\node (v1) at (0,0) {$S$};
\node (v2) at (3,0) {$S$};
\node (v3) at (6,0) {$S$};
\node (v4) at (9,0) {$S$};
\node (v5) at (12,0) {$S$};

\draw (v1) -- (v2);
\draw (v2) -- (v3);
\draw (v3) -- (v4);
\draw (v1) to[bend right =15](v3);
\draw (v2) to [bend right = 15] (v4);
\draw (v1) to [bend right = 20] (v4);
\draw (v1) to [bend right = 25] (v5);
\draw (v2) to [bend right = 20] (v5);
\draw (v3) to [bend right = 15] (v5);
\draw (v4) to (v5);

\node (u1) at (0,3) {$K$};
\node (u2) at (3,3) {$K$};
\node (u3) at (6,3) {$K$};
\node (u4) at (9,3) {$K$};
\node (u5) at (12,3) {$K$};

\draw (u1) -- (v1);
\draw (u2) -- (v2);
\draw (u3) -- (v3);
\draw (u4) -- (v4);
\draw (u5) -- (v5);

\tikzstyle{every node}=[above right,outer sep=3pt]
\draw (v1) node []           {$Z_1$};
\draw (v2) node []           {$Z_2$};
\draw (v3) node []           {$Z_3$};
\draw (v4) node []           {$Z_4$};
\draw (v5) node []           {$Z_5$};
\tikzstyle{every node}=[right,outer sep=6pt]
\draw (u1) node []           {$Y_1$};
\draw (u2) node []           {$Y_2$};
\draw (u3) node []           {$Y_3$};
\draw (u4) node []           {$Y_4$};
\draw (u5) node []           {$Y_5$};

}

\end{tikzpicture}

\caption{A candelabrum. (Lines indicate complete pairs, and $K,S$ mark cliques and stable sets.)} \label{fig:candelabra}
\end{figure}

Let $G$ be a graph, let $k\ge 1$ be an integer, and let $Y_1\ll Y_k, Z_1\ll Z_k$ be disjoint nonempty sets of $V(G)$
with union $V(G)$, such that
\begin{itemize}
\item $Y_1\ll Y_k$ are cliques, and $Z_1\ll Z_k$ are stable sets;
\item for $1\le i<j\le k$,  $Y_i$ is anticomplete to $Y_j$, and $Z_i$ is complete to $Z_j$;
\item for $1\le i,j\le k$, $Y_i$ is complete to $Z_j$ if $i=j$, and otherwise $Y_i$ is anticomplete to $Z_j$.
\end{itemize}
We call such a graph $G$ a {\em candelabrum}, with {\em base} $Z_1\cup\cdots\cup Z_k$. (See figure \ref{fig:candelabra}.)

Candelabra are uncluttered, but they can also be used to make larger uncluttered graphs. Let $H_1$ be an uncluttered graph,
and let $H_2$ be a candelabrum with base $Z$. Take the disjoint union of $H_1,H_2$, and add edges to make $V(H_1)$ complete
to $Z$. Let $G$ be the graph we produce. Then $G$ is uncluttered (we leave checking this to the reader), 
and $V(H_1)$ is a homogeneous set of $G$. We say $G$ is {\em candled} if it can be constructed by this process; that is,
if $G$ has an induced subgraph $H$ that is a candelabrum, with base $Z$ say, and $Z$ is complete to $V(G)\setminus V(H)$,
and $V(H)\setminus Z$ is anticomplete to $V(G)\setminus V(H)$. (Thus $V(G)\setminus V(H)$ is a homogeneous set, but it might
not be nontrivial; it might even be empty.)

Now we can state our main theorem:
\begin{thm}\label{mainthm}
Let $G$ be an uncluttered graph. Then either
\begin{itemize}
\item one of $G,\overline{G}$ is disconnected; or
\item one of $G, \overline{G}$ has adjacent simplicial twins; or
\item one of $G,\overline{G}$ is candled; or
\item one of $G,\overline{G}$ is the line graph of a triangle-free graph.
\end{itemize}
\end{thm}
The proof of \ref{mainthm} will occupy the remainder of the paper, and will be completed in \ref{mainthm2}.
\section{Homogeneous sets}

In this section we will prove:
\begin{thm}\label{nohomog}
Let $G$ be an uncluttered graph. Suppose that
\begin{itemize}
\item $G$ is connected and anticonnected;
\item $G$ has no adjacent simplicial twins, and no nonadjacent antisimplicial twins; and
\item $G,\overline{G}$ are not candled.
\end{itemize}
Then $G$ has no nontrivial homogeneous set.
\end{thm}
The proof requires several steps, that we carry out in this section. 
Let us say
an {\em anticomponent} of a graph $G$ is a induced subgraph whose complement graph is a component of $\overline{G}$.
We begin with:
\begin{thm}\label{mixedhom}
Let $G$ be uncluttered, and connected and anticonnected, with a nontrivial homogeneous set 
that is not a clique or stable set. Then one of $G,\overline{G}$ is candled.
\end{thm}
\Proof
Let $X$ be a nontrivial homogeneous set that is not a clique or stable set. 
By replacing $X$ by a superset if necessary, we may assume that no proper superset of $X$ is a nontrivial homogeneous set.
Let $Z$ be the set of vertices in $V(G)\setminus X$ that are complete to $X$, and let $Y$ be the set of vertices in $V(G)\setminus X$
that are anticomplete to $X$. Let $Y_1\ll Y_k$ be the vertex sets of the components of $G[Y]$, and let $Z_1\ll Z_{\ell}$ be the 
vertex sets of the anticomponents of $G[Z]$.
\\
\\
(1) {\em $Y_1\ll Y_k$ and $Z_1\ll Z_{\ell}$ are homogeneous sets of $G$.}
\\
\\
Suppose that $Y_1$ is not a homogeneous set; then there exists $v\in V(G)\setminus Y_1$ that has a neighbour and a non-neighbour in $Y_1$,
and since $G[Y_1]$ is connected, there is an edge $yy'$ of $G[Y_1]$ such that $v$ is adjacent to $y$ and not to $y'$. Since $X$
is anticomplete to $Y_1$, it follows that $v\notin X$, and similarly $v\notin Y_2\ll Y_k$; so $v\in Z$. Now $X$ is not a clique,
so there exist nonadjacent $x_1,x_2\in X$. But then $G[\{x_1,x_2,v,y,y'\}]$ is a fork, a contradiction. Thus
$Y_1$ is a homogeneous set, and similarly so are $Y_2\ll Y_k$; and so are $Z_1\ll Z_{\ell}$, by applying the same argument in the 
complement graph. This proves (1).

\bigskip
It follows that for each $Y_i$ and each $Z_j$, $Y_i$ is either complete or anticomplete to $Z_j$.
\\
\\
(2) {\em $Y,Z$ are both nonempty. Moreover, for $1\le i\le k$, there exist $j,j'\in \{1\ll \ell\}$ such that $Y_i$ is complete to $Z_j$
and anticomplete to $Z_{j'}$, and for $1\le j\le \ell$, there exist $i,i'\in \{1\ll k\}$ such that $Z_{j}$ is complete to $Y_i$
and anticomplete to $Y_{i'}$.}
\\
\\
Since $X$ is a nontrivial homogeneous set, $X\ne V(G)$ and so $Y\cup Z\ne \emptyset$. Since $G$ is connected, it follows that
$Y\ne \emptyset$,
and similarly $Z\ne \emptyset$ since $G$ is anticonnected. Since $G$ is connected, for $1\le i\le k$ there exists
$j\in \{1\ll \ell\}$ such that $Y_i$ is complete to $Z_j$; and also, since $X\cup \{Y_i\}$ is not a homogeneous set from the maximality
of $X$, it follows that there exists $j'\in \{1\ll \ell\}$ such that $Y_i$ is anticomplete to $Z_{j'}$. 
The same argument in the complement shows the final statement. This proves (2).

\bigskip
Let $i_1\ll i_t\in \{1\ll k\}$ be distinct, and let $j_1\ll j_t\in \{1\ll \ell\}$ be distinct. We say the pairs $(i_1,j_1)\ll (i_t,j_t)$
form a {\em matching of order $t$} if for $1\le r,s\le t$, $Y_{i_r}$ is complete to $Z_{j_s}$ if $r=s$ and otherwise $Y_{i_r}$ is 
anticomplete to $Z_{j_s}$. Similarly the pairs form an {\em antimatching of order $t$} if for $1\le r,s\le t$, $Y_{i_r}$ is 
anticomplete to $Z_{j_s}$ if $r=s$ and otherwise $Y_{i_r}$ is
complete to $Z_{j_s}$.
\\
\\
(3) {\em There exist pairs $(i,j), (i',j')$  that form a matching of order two.}
\\
\\
Choose $i\in \{1\ll k\}$ such that $Y_i$
is complete to $Z_j$ for as many $j\in \{1\ll \ell\}$ as possible; by (2), we can choose $j'\in \{1\ll \ell\}$ 
such that $Y_i$ is anticomplete to
$Z_{j'}$; by (2) we can choose $i'\in \{1\ll k\}$ such that $Y_{i'}$ is complete to $Z_{j'}$; and then 
from the choice of $i$, there exists $j\in \{1\ll\ell\}$
such that $Z_j$ is complete to $Y_i$ and not to $Y_{i'}$. Then $(i,j),(i',j')$ forms a matching of order two. This proves (3).

\bigskip

Choose $t\ge 2$ maximum such that there are $t$ pairs $(i,j)$ (with $i\in \{1\ll k\}$ and $j\in \{1\ll \ell\}$) 
that form a matching or antimatching; and by 
taking complements if necessary, we may assume the pairs form a matching. By renumbering, we may assume that $(1,1)\ll (t,t)$
form a matching. For $1\le i\le t$, choose $y_i\in Y_i$ and $z_i\in Z_i$.
\\
\\
(4) {\em Every vertex in $Y_{t+1}\cup \cdots\cup Y_k$ is complete or anticomplete to $Z_1\cup\cdots\cup Z_t$, and 
every vertex in $Z_{t+1}\cup \cdots\cup Z_{\ell}$ is complete or anticomplete to $Y_1\cup\cdots\cup Y_t$.}
\\
\\
If $v\in Y_{t+1}\cup \cdots\cup Y_k$ has a neighbour in $Z_1$ and a nonneighbour in $Z_2$ say, then 
$G[\{v,z_1,z_2,y_1,y_2\}]$ is a fork, a contradiction. Similarly,
if $v\in Z_{t+1}\cup \cdots\cup Z_{\ell}$ has a neighbour in $Y_1$ and a nonneighbour in $Y_2$, 
then $G[\{v,z_1,z_2,y_1,y_2\}]$ is an antifork, a contradiction.
This proves (4).

\bigskip

Let $P$ be the set of $v\in Y_{t+1}\cup \cdots\cup Y_k$ such that $v$ is complete to $Z_1\cup\cdots\cup Z_t$, and let
$P'$ be the set of $v\in Y_{t+1}\cup \cdots\cup Y_k$ such that $v$ is anticomplete to $Z_1\cup\cdots\cup Z_t$.
Let $Q$ be the set of $v\in Z_{t+1}\cup \cdots\cup Z_{\ell}$ such that $v$ is complete to $Y_1\cup\cdots\cup Y_t$, and let
$Q'$ be the set of $v\in Z_{t+1}\cup \cdots\cup Z_{\ell}$ such that $v$ is anticomplete to $Y_1\cup\cdots\cup Y_t$.
\\
\\
(5) {\em $X\cup Y_1\cup\cdots\cup Y_t\cup Z_1\cup\cdots\cup Z_t\cup P\cup Q'$ is a homogeneous set.}
\\
\\
Let us call this set $A$. We will show that $P'$ is anticomplete to $A$, and $Q$ is complete to $A$. Let $v\in P'$; then
$v$ is anticomplete to $X$ since $v\in Y$; $v$ is anticomplete to $Y_1\cup\cdots\cup Y_t$ since $Y_1\ll Y_t$
are vertex sets of components of $G[Y]$; and $v$ is anticomplete to $Z_1\cup\cdots\cup Z_t$ by definition of $P'$. Since
$v$ belongs to some $Y_i$, and $Y_i$ is a homogeneous set by (1), and so $Y_i\subseteq P'$, it follows that $v$
is anticomplete to $P$. Also from the maximality of $t$, $v$ is anticomplete to $Q'$. This proves that $P'$
is anticomplete to $A$.

Now let $v\in Q$. We must show that $v$ is complete to $A$. Certainly $v$ is complete to $X$ and to $Z_1\cup\cdots\cup Z_t$, 
and to $Y_1\cup\cdots\cup Y_t$ from the definition of $Q$.
Since $v$ belongs to an anticomponent of $G[Z]$ with vertex set in $Q$, it follows that $v$ is complete to $Q'$. It remains to show
that $v$ is complete to $P$. Suppose not, and let $u\in P$ be nonadjacent to $v$. For $1\le i\le t$, choose $y_i\in Y_i$ 
and choose $z_i\in Z_i$. If $t\ge 3$, $G[\{v,y_1,y_2,u,z_3\}]$
is a fork, a contradiction, and so $t=2$. Let $u\in Y_3$ and $v\in Z_3$ say; then $Y_3$ is anticomplete to $Z_3$ by (1), and so
the three pairs $(1,2),(2,1),(3,3)$ form an antimatching, contrary to the maximality of $t$. This proves (5).

\bigskip
From the maximality of $X$, (5) implies that the set of (5) is not a nontrivial homogeneous set, and so $P', Q=\emptyset$.
But then $X\cup P\cup Q'$ is a nontrivial homogeneous set, and the maximality of $X$ implies that $P,Q'=\emptyset$. 
If there exist nonadjacent $y,y'\in Y_1$, then $G[\{y,y',z_1,z_2,y_2\}]$ is a fork, a contradiction; so $Y_1\ll Y_t$ are cliques.
If there exist adjacent $z,z'\in Z_1$, then $G[\{z,z',y_1,z_2,y_2\}]$ is an antifork, a contradiction. Thus $Z_1\ll Z_t$
are stable sets.
But then $G$ is candled, as required. This proves \ref{mixedhom}.~\bbox

Let $P$ be a four-vertex induced path in a graph $G$. A {\em centre} for $P$ means a vertex of $V(G)\setminus V(P)$ that is complete
to $V(P)$, and an {\em anticentre} for $P$ is a vertex of $V(G)\setminus V(P)$ that is anticomplete
to $V(P)$.
\begin{thm}\label{centre+anti}
Let $G$ be uncluttered, and let $P$ be a four-vertex induced path in $G$. If there is a centre and an anticentre for $P$ then
there is a nontrivial homogeneous set in $G$ that is not a clique or stable set.
\end{thm}
\Proof
Let $A$ be the set of all anticentres for $P$, and let $C$ be the set of all centres for $P$. Thus $A,C\ne \emptyset$. Let $B$
be the set of $v\in V(G)\setminus (V(P)\cup A\cup C)$; thus $B$ is the set of all vertices not in $V(P)$ with a neighbour and a 
nonneighbour in $V(P)$. Now either every vertex in $C$ has a neighbour in $A$, or every vertex in $A$ has a nonneighbour
in $C$; and by taking complements if necessary, we assume the first. Let $P$ have vertices $p_1\dd p_2\dd p_3\dd p_4$ in order.
\\
\\
(1) {\em $B$ is complete to $C$.}
\\
\\
Let $b\in B$ and $c\in C$, and suppose that $b,c$ are nonadjacent. Choose $a\in A$ adjacent to $c$. Suppose first that $a,b$ are 
nonadjacent. Since $P$ is anticonnected, and $b\in B$, there exist nonadjacent $p,p'\in V(P)$ such that $b$ is adjacent to $p$
and not to $p'$; but then $G[\{p,p',a,b,c\}]$ is a fork, a contradiction. So $a,b$ are adjacent. Let $I$
be the set of $i\in \{1\ll 4\}$ such that $b,p_i$ are adjacent. Since $G[\{a,b,c,p_1,p_4\}]$
is not a fork, one of $1,4\in I$, and we may assume $1\in I$ without loss of generality. If $4\notin I$, then $2\in I$ since
$G[\{a,b,c,p_2,p_4\}]$ is not a fork; so $3\notin I$ since the subgraph induced on $G[\{b,p_1,p_2,p_3, p_4\}]$ 
is not an antifork; but then $G[\{b,c,p_1,p_2,p_4\}]$ is an antifork, a contradiction. Thus $4\in I$. Since $b$
is not a centre, one of $2,3\notin I$, and we assume $3\notin I$ without loss of generality. But then
$G[\{a,b,p_1,p_3,p_4\}]$ is a fork, a contradiction. This proves (1).

\bigskip
Let $A'$ be the union of the vertex sets of the components of $G[A]$ that are not anticomplete to $B$, and let $A''=A\setminus A'$.
\\
\\
(2) {\em $A'$ is complete to $C$.}
\\
\\
Let $c\in C$, and suppose that $c$ is not complete to $A'$. From the definition of $A'$, there is an induced path $Q$
with one end in $B$, and all other vertices in $A'$, such that some vertex of $Q$ is not adjacent to $c$. Choose $Q$ minimal, with
vertices $q_1\cc q_k$ in order, where $q_k\in B$ and $q_1\ll q_{k-1}\in A$. From (1), $k\ge 2$.
From the minimality of $Q$, it follows that $c$
is nonadjacent to $q_1$ and adjacent to all of $q_2\ll q_{k}$. Choose adjacent $q_{k+1}, q_{k+2}$ of $P$ such that $q_k$ is adjacent
to $q_{k+1}$ and not to $q_{k+2}$. (This is possible since $q_k\in B$.) Then $G[\{c,q_1,q_2,q_3,q_4\}]$ is an antifork, a contradiction.
This proves (2).

\bigskip

Let $X=V(P)\cup B\cup A'$. From (1) and (2), $C$ is complete to $X$; and from the definition of $A''$, $A''$ is anticomplete to $X$.
Thus $X$ is a nontrivial homogeneous set satisfying the theorem. This proves \ref{centre+anti}.~\bbox

\begin{thm}\label{newlemma}
Let $G$ be uncluttered and connected, and let $A,B,C,D$ be disjoint subsets of $V(G)$, with union $V(G)$, and with the following properties:
\begin{itemize}
\item $A$ is a clique and $A\ne \emptyset$;
\item $B$ is a stable set and $|B|\ge 2$;
\item $A$ is complete to $B$ and anticomplete to $C,D$;
\item $B$ is complete to $C$ and anticomplete to $D$.
\end{itemize}
Then $G$ is candled.
\end{thm}
\Proof
Let $C_1\ll C_k$ be the vertex sets of the anticomponents of $G[C]$, and let $D_1\ll D_{\ell}$ be the vertex sets of the 
components of $G[D]$.
\\
\\
(1) {\em $C_1\ll C_k$ and $D_1\ll D_{\ell}$ are homogeneous sets.}
\\
\\
Suppose $C_1$ is not a homogeneous set; then there exist $d\in V(G)\setminus C_1$ 
and nonadjacent $c,c'\in C_1$ with $d$ adjacent to $c$ and not to $c'$.
Choose $a\in A$ and $b\in B$; then $G[\{a,b,c,c',d\}]$ is a fork, a contradiction. Thus $C_1\ll C_k$ are all homogeneous sets.

Now suppose $D_1$ is not a homogeneous set; 
then similarly there exists $c\in C$ and adjacent $d,d'\in D_1$ such that $c$ is adjacent to $d$ and not to $d'$. Since
$|B|\ge 2$ and $B$ is stable, there exist nonadjacent $b,b'\in B$; but then $G[\{b,b',c,d,d'\}]$ is a fork, a contradiction.
This proves (1).

\bigskip
It follows that for each $C_i$ and each $D_j$, $C_i$ is either complete or anticomplete to $D_j$.
For $1\le i\le k$ choose $c_i\in C_i$, and for $1\le j\le \ell$ choose $d_j\in D_j$. Choose $a\in A$ and $b\in B$.
\\
\\
(2) {\em For $1\le j\le \ell$, there is a unique value of $i\in \{1\ll k\}$ such that $D_j$ is complete to $C_i$.}
\\
\\
We assume $j=1$ without loss of generality. Since $G$ is connected and $A\cup B\cup D_2\cup\cdots\cup D_{\ell}$ is anticomplete
to $D_1$, it follows from (1) that there exists $i\in \{1\ll k\}$ such that $D_1$ is complete to $C_i$. 
Suppose there are two such values
of $i$, say $i=1,2$; then $G[\{a,b,c_1,c_2,d_1\}$ is an antifork, a contradiction. This proves (2).
\\
\\
(3) {\em For $1\le i\le k$ there is at most one value of $j\in \{1\ll \ell\}$ such that $C_i$ is complete to $D_j$.}
\\
\\
For let $i = 1$ say, and suppose that $C_1$ is complete to $D_1,D_2$ say. Then $G[\{a,b,c_1,d_1,d_2\}]$ is a fork, a contradiction.
This proves (3).

\bigskip

From (2) and (3) it follows that $k\ge \ell$, and we may renumber such that $D_i$ is complete to $C_i$ for $1\le i\le \ell$.
For $1\le i\le k$, $D_i$ is a clique, since if $d,d'\in D_i$ are nonadjacent then $G[\{a,b,c_i,d,d'\}]$ is a fork; and
also $C_i$ is a stable set, since if $c,c'\in C_i$ are adjacent then $G[\{a,b,c,c',d_i\}]$ is an antifork. 
Thus the restriction of $G$ to 
$$A\cup B\cup C_1\cup\cdots\cup C_{\ell}\cup D_1\cup\cdots\cup D_{\ell}$$
is a candelabrum with base $B\cup C_1\cup\cdots\cup C_{\ell}$.
Since $C_{\ell+1}\cup\cdots\cup C_k$ is complete to $B\cup C_1\cup\cdots\cup C_{\ell}$ and anticomplete to
$A\cup D_1\cup\cdots\cup D_{\ell}$, it follows $G$ is candled.
This proves \ref{newlemma}.~\bbox

Now we can prove the main result of this section, which we restate:
\begin{thm}\label{nohomog2}
Let $G$ be an uncluttered graph. Suppose that
\begin{itemize}
\item $G$ is connected and anticonnected;
\item $G$ has no adjacent simplicial twins, and no nonadjacent antisimplicial twins; 
\item $G,\overline{G}$ are not candled.
\end{itemize}
Then $G$ has no nontrivial homogeneous set.
\end{thm}
\Proof For each nontrivial homogeneous set $X$, we define its ``score'' as follows. By 
\ref{mixedhom}, $X$ is either a clique or a stable set. If $X$ is stable, its {\em score} is the number of components of
$G[Y]$, where $Y$ is the set of all vertices in $V(G)\setminus X$ that are anticomplete to $X$. If $X$ is a clique, its
{\em score} is the number of anticomponents of $G[Z]$, where $Z$ is the set of vertices in $V(G)\setminus X$ that are 
complete to $X$.

Suppose that there is a nontrivial homogeneous set $X$, and choose $X$ with minimum score. 
Let $Z$ be the set of vertices in $V(G)\setminus X$ that are complete to $X$,
and $Y$ the set that are anticomplete to $X$.
Let $Y_1\ll Y_k$ be the vertex sets of the components of $G[Y]$, and let 
$Z_1\ll Z_{\ell}$ be the vertex sets of the anticomponents of $G[Z]$. By replacing $G$ by its complement if necessary, we may assume that
$X$ is a stable set.
\\
\\
(1) {\em $Y_1\ll Y_k$ are homogeneous sets, and cliques.}
\\
\\
If $Y_1$ is not a homogeneous set, then there exists $z\in Z$ and $y,y'\in Y_1$, such that $z\dd y\dd y'$ is an induced path (because
$G[Y_1]$ is connected). Let $x,x'\in X$ be nonadjacent; then $G[\{x,x',z,y,y'\}]$ is a fork, a contradiction. This proves that
$Y_1$ is a homogeneous set. By \ref{mixedhom}, $Y_1$ is either a stable set or a clique; but $G[Y_1]$ is connected, and so
$Y_1$ is a clique. This proves (1).

\bigskip
By hypothesis, there are no nonadjacent antisimplicial twins in $G$; and in particular, vertices in $X$ are not antisimplicial.
Thus $Y$ is not stable, and so we may assume that $|Y_1|\ge 2$. Let $N$ be the set of vertices in $Z$ that are complete
to $Y_1$. By hypothesis, there are no adjacent simplicial twins, and in particular the vertices in $Y_1$ are not simplicial;
so $N$ is not a clique. Hence there is an anticomponent $D$ of $N$ with at least two vertices. Let $N'=N\setminus D$, and
$Y'=Y_2\cup\cdots\cup Y_k$. Thus the six sets $X,D, N', Z\setminus N, Y_1, Y'$  are pairwise disjoint and have union $V(G)$.
\\
\\
(2) {\em $D$ is a homogeneous set, and $D$ is stable.}
\\
\\
Since $G[D]$ is anticonnected, if $D$ is not a homogeneous set then there exists $v\in V(G)\setminus D$ and nonadjacent $d,d'\in D$,
such that $v$ is adjacent to $d$ and not to $d'$. Since $X\cup N'\cup Y_1$ is complete to $D$, it follows that $v\in (Z\setminus N)\cup Y'$.
Choose $y,y'\in Y_1$, adjacent; then 
$G[\{v,d,d',y,y'\}]$ is an antifork, a contradiction. Thus $D$ is homogeneous. By \ref{mixedhom} is either stable or a clique,
and it is not a clique since $G[D]$ is anticonnected and $|D|\ge 2$. This proves (2).

\bigskip

Let $U$ be the set of vertices in $Z\setminus N$ that are complete to $D$, and let $W$ be the set of vertices in $Z\setminus N$
that are anticomplete to $D$. Thus $U\cup W=Z\setminus N$ by (2).
\\
\\
(3) {\em $W$ is a clique.}
\\
\\
If $w,w'\in W$ are nonadjacent, choose $y\in Y_1$, $d\in D$ and $x\in X$; then $G[\{w,w',y_1,d_,x\}]$ is a fork, a contradiction.
This proves (3).
\\
\\
(4) {\em $Y'$ is anticomplete to $D\cup W$, and $N'$ is complete to $W$, and $U$ is anticomplete to $W$.}
\\
\\
From the minimality of the score of $X$, it follows that the score of $D$ is at least that of $X$, that is, at least $k$.
For $2\le i\le k$, $Y_i$ is a homogeneous set by (1), and since $D$ is also a homogeneous set by (2), it follows that $Y_i$
is complete or anticomplete to $D$. Let $C$ denote the set of vertices not in $D$ and anticomplete to $D$. Then $C$
consists of $W$ and some of $Y_2\ll Y_k$, those $Y_i$ that are anticomplete to $D$. In particular, since $W$ is a clique
by (3), 
the score of $D$ is one more (for $W$) than the number of $i\in \{2\ll k\}$
such that $Y_i$ is anticomplete to $D\cup W$. Since the score of $D$ is at least $k$, it follows that $Y_2\ll Y_k$ are all anticomplete
to $D\cup W$. This proves the first assertion of (4). 

For the second assertion of (4), let $y\in Y_1, d\in D, x\in X$ and $w\in W$; then for $n\in N'$, if $n$ is nonadjacent to $w$
then $G[\{y,d,x,w,n\}]$ is an antifork. This proves that $N'$ is complete to $W$. 

For the third assertion, let $u\in U$, and let $y,d,x,w$ be as before; if $u,w$ are adjacent then $G[\{y,d,x,u,w\}]$
is an antifork, a contradiction. This proves (4).

\bigskip

If $N'=\emptyset$, then the four sets $W,X,D\cup U, Y$ (in this order) satisfy the hypotheses of \ref{newlemma}, a contradiction.
Thus $N'\ne \emptyset$. If $Y'=\emptyset$, then the four sets $X,Y_1,W\cup U, N'\cup D$ satisfy the hypotheses of \ref{newlemma}
applied in $\overline{G}$, a contradiction. So $Y'\ne \emptyset$. Choose $y\in Y_1$, $y'\in Y'$, $d\in D$, $x\in X$,
$n\in N'$  and $w\in W$; then $y\dd d\dd x\dd w$ is a four-vertex induced path, and $n$ is a centre and $y'$ is an anticentre,
contrary to \ref{centre+anti}. This proves \ref{nohomog2}.~\bbox

\section{Graphs without homogeneous sets}

In view of \ref{nohomog}, henceforth we can restrict our attention to uncluttered graphs with no nontrivial homogeneous set.
First here is a useful lemma.
\begin{thm}\label{mixedvert}
Let $G$ be a graph with no nontrivial homogeneous set, and let $A\subseteq V(G)$, not a clique, with $A\ne V(G)$. Then there 
exist $v\in V(G)\setminus A$ and nonadjacent $a,a'\in A$, such that $v$ is adjacent to $a$ and not to $a'$.
\end{thm}
\Proof
Since $A$ is not a clique, there is an anticomponent $X$ of $G[A]$ with at least two vertices. Since $X\ne V(G)$
and $X$ is not a nontrivial homogeneous set, there is a vertex $v\in V(G)\setminus X$ with a neighbour and a nonneighbour in $X$.
Since $X$ is an anticomponent of $G[A]$, it follows that $v\notin A$. Since $G[X]$ is anticonnected, and $v$ has a neighbour
and a nonneighbour in $X$, there exist
nonadjacent $a,a'\in X$ such that $v$ is adjacent to $a$ and not to $a'$. This proves \ref{mixedvert}.~\bbox

If $X$ is a subgraph of $G$, or a set of vertices of $G$, we say that $X$ is {\em dominating} if every vertex not in $X$
has a neighbour in $X$. The {\em diamond} is the graph with four vertices and five edges, and a {\em triangle}
is a clique of cardinality three. We begin with:

\begin{thm}\label{diamond}
Let $G$ be an uncluttered graph with no nontrivial homogeneous set. Then every diamond in $G$ is dominating. Moreover,
every triangle that is contained in a diamond of $G$ is dominating.
\end{thm}
\Proof
Suppose that there is a diamond $D$ that is not dominating. Let $v_1,v_2$ be the two vertices of $D$ that have degree three
in $D$; so there is a set $A$ of $G$, not a clique, such that $v_1,v_2\notin A$ and $A$ is complete to $\{v_1,v_2\}$,
and $A\cup \{v_1,v_2\}$ is not dominating. Choose such a set $A$, maximal. Since $A$ is not a nontrivial homogeneous set,
\ref{mixedvert} implies that there exist $v\in V(G)\setminus A$ and nonadjacent $a,a'\in A$
such that $v$ is adjacent to $a$ and not to $a'$. Thus $v\ne v_1,v_2$.
If $v$ is nonadjacent to both $v_1,v_2$, then $G[\{v,v_1,v_2,a,a'\}]$ is an antifork, a contradiction; so we may
assume that $v$ is adjacent to $v_1$. Since $A\cup \{v_1,v_2\}$ is not dominating, there is a vertex $w$ that is anticomplete
to $A\cup \{v_1,v_2\}$. Now there are four possibilities: $w$ may or may not be adjacent to $v$, and $v$ may or may not
be adjacent to $v_2$. Suppose first that $v$ is not adjacent to $v_2$, and so $v\dd a\dd v_2\dd a'$ is an induced four-vertex
path, $P$ say. Now $v_1$ is a centre for $P$, so by \ref{centre+anti}, $w$ is not an anticentre, and therefore $w$ is adjacent to $v$.
But then $G[\{w,v,v_1,a,v_2\}]$ is an antifork, a contradiction. So $v$ is adjacent to $v_2$. If $w,v$ are nonadjacent,
then we can add $v$ to $A$, contrary to the maximality of $A$. Thus $v,w$ are adjacent; but then 
$G[\{w,v,v_1,v_2,a'\}]$ is an antifork,
a contradiction. This proves that every diamond is dominating.

Now let $T$ be a triangle, contained in a diamond; thus some vertex $v\in V(G)\setminus T$ has exactly two neighbours
in $T$. Suppose that $T$ is not dominating, and let $w$ be anticomplete to $T$. Since the diamond $G[T\cup \{v\}]$
is dominating, $w$ is adjacent to $v$; but then $G[T\cup \{v,w\}]$ is an antifork, a contradiction. This 
proves \ref{diamond}.~\bbox

\begin{thm}\label{nondomtri}
Let $G$ be an uncluttered graph with no nontrivial homogeneous set,
such that $G$ is not the line graph of a bipartite graph.
Then for every nondominating triangle $T$ in $G$, there is a unique maximal clique including $T$, say $C$, and it is not dominating.
Moreover, for each $v\in C$, if $C$ is not the only maximal clique containing $v$, then there is exactly one other
maximal clique containing $v$, say $C_v$, and $C_v\cap C=\{v\}$, and $C_v$ is not dominating.
\end{thm}
\Proof
Let $C$ be a maximal nondominating clique including $T$; and let $S$ be the set of vertices that are anticomplete to $C$.
Thus $S\ne \emptyset$. Let $A$ be the set of vertices in $V(G)\setminus C$ with exactly one neighbour in $C$, and let $B$
be the set of vertices in $V(G)\setminus C$ with at least two neighbours in $C$. Thus $A,B,C,S$ are pairwise disjoint and
have union $V(G)$. 
\\
\\
(1) {\em $B$ is complete to $C\cup S$.}
\\
\\
No triangle included in $C$ is dominating, and hence by \ref{diamond}, no triangle included in $C$
is contained in a diamond. Consequently every vertex in $B$ is complete to $C$. From the maximality of $C$, 
$B$ is complete to $S$. This proves (1).

\bigskip

Every vertex in $A$ has a unique neighbour in $C$. For each $c\in C$, let $A_c$ be the set of vertices in $A$ that are
adjacent to $c$. Thus $A=\bigcup_{c\in C}A_c$. 
\\
\\
(2) {\em There is at most one $c\in C$ with $A_c=\emptyset$.}
\\
\\
For otherwise the set of $c\in C$ with $A_c=\emptyset$ is a nontrivial homogeneous set, which is impossible. This proves (2).
\\
\\
(3) {\em $|B|\le 1$, and $B$ is anticomplete to $A$.}
\\
\\
Suppose first that $B$ is not a clique. By \ref{mixedvert},
there is a vertex $v\in V(G)\setminus B$ and nonadjacent $b,b'\in B'$ such that $v$ is adjacent to 
$b$ and not to $b'$. By (1), $v\in A_c$
for some $c\in C$. Choose $c_1,c_2\in C$ different from $c$; then $G[\{v,b,b',c_1,c_2\}]$ is an antifork, a contradiction.
This proves that $B$ is a clique.

Suppose that some $b\in B$ has a neighbour in $A$.
Since $G$ is anticonnected (because it has no nontrivial homogeneous set), it follows that $b$ has a nonneighbour;
and since $b$ is complete to $C\cup S$ and $B$ is a clique, it follows that $b$ has a nonneighbour in $A$. By (2), there are at least
two vertices $c\in C$ such that $A_c\ne \emptyset$; and so there exist distinct $c_1,c_2\in C$, and $a_i\in A_{c_i}$
for $i = 1,2$, such that $b$ is adjacent to $a_1$ and nonadjacent to $a_2$. (To see this, choose a neighbour $a\in A$ of $b$
and a nonneighbour $a'\in A$ of $b$. If $a,a'$ both belong to $A_{c_1}$ say, choose $a''\in A_{c_2}$, and replace 
one of $a,a'$ by $a''$.) Choose $c_3\in C\setminus \{c_1,c_2\}$.
Since $G[\{b,c_1,c_3,a_1\}]$ is a diamond, the triangle $G[\{b,c_1,c_3\}]$
is dominating by \ref{diamond}, and yet $a_2$ has no neighbour in this triangle, a contradiction.

Thus $B$ is anticomplete to $A$, and so $B$ is a homogeneous set; and hence $|B|\le 1$.
This proves (3).
\\
\\
(4) {\em For each $c\in C$, $A_c$ is a clique; and if $B\ne \emptyset$ then $S$ is a clique.}
\\
\\
Suppose that $c\in C$ and $A_c$ is not a clique; then by \ref{mixedvert}, 
there exists $v\notin A_c$ and nonadjacent $a,a'\in A_c$ such that
$v$ is adjacent to $a$ and not to $a'$. Thus $v\in A\cup S$, by (3); and so there exists $c'\in C\setminus \{c\}$ nonadjacent
to $v$. But then $G[\{c,c',a,a',v\}]$ is a fork, a contradiction. Thus each $A_c$ is a clique. Now suppose that $B=\{b\}$ say,
and $S$ is not a clique. Then by \ref{mixedvert}, there exists $v\notin S$ and nonadjacent $s,s'\in S$ such that $v$ is
adjacent to $s$ and not to $s'$. Thus $v\in A_c$ for some $c\in C$. Choose $c'\in C\setminus \{c\}$. Then $G[\{b,s,s',v,c'\}]$
is a fork, a contradiction. This proves (4).
\\
\\
(5) {\em $B= \emptyset$.}
\\
\\
Suppose that $B\ne \emptyset$, and so $|B|=1$ by (3). Let $B=\{b\}$ say.
Define $A_b=S$; thus $C\cup \{b\}$ is a clique $D$ say, and every vertex not in $D$ belongs to one
of the sets $A_d\;(d\in D)$, and has
a unique neighbour in $D$; and each set $A_d\;(d\in D)$ is a clique by (4). We claim:
\begin{itemize}
\item for all distinct $d_1,d_2\in D$, each vertex in $A_{d_1}$ has at most one neighbour in $A_{d_2}$; and
\item for all distinct $d_1,d_2,d_3\in D$, if $a_i\in A_{d_i}$ for $i = 1,2,3$, and $a_1$ is adjacent
to both $a_2,a_3$, then $a_2,a_3$ are adjacent.
\end{itemize}
To see the first claim, 
suppose that $v\in A_{d_1}$ has two neighbours $a,a'\in A_{d_2}$. By (4), $a,a'$ are adjacent. 
Choose $d_3\in D\setminus \{d_1,d_2\}$; then 
$G[\{v,a,a',d_2,d_3\}]$ is an antifork, a contradiction. This proves the first claim.

For the second claim, 
suppose that $d_1,d_2,d_3\in D$ are distinct, and $a_i\in A_{d_i}$ for $i = 1,2,3$, and $a_1$ is adjacent to $a_2,a_3$,
and $a_2,a_3$ are not adjacent. Since $|C|\ge 3$ it follows that $|D|\ge 4$; choose $d_4\in D\setminus \{d_1,d_2,d_3\}$. Then 
$G[\{a_1,a_2,a_3,d_1,d_4\}]$ is a fork, a contradiction. This proves the second claim.

Let $H$ be the subgraph of $G$ obtained by deleting the edges of the cliques
$A_d\cup \{d\}\;(d\in D)$. From the two bullet claims above, it follows that each component of $H$ is a clique. 
($D$ itself is one such component.)
Thus we have found two sets of cliques of $G$; the sets $A_d\cup \{d\}\;(d\in D)$, and the components of $H$. Each vertex
of $G$ belongs to exactly one clique in the first set, and exactly one in the second; and every edge of $G$
belongs to one of the cliques in one of the sets. Consequently $G$ is the line graph
of a bipartite graph, a contradiction. This proves that $B=\emptyset$, and so proves (5).

\bigskip
It follows that $C$ is a maximal clique of $G$, and is nondominating. Moreover, every edge of $C$ is not contained in any other maximal
clique; and every vertex $c\in C$ is contained in at most two maximal cliques, namely $C$ and $A_c\cup\{c\}$ if $A_c\ne \emptyset$.
To complete the proof of the theorem, we only need to show that the cliques $A_c\cup \{c\}$ are not dominating. Suppose then
that $c\in C$, and $A_c\cup \{c\}$ is dominating. We have already shown that every clique including a nondominating triangle is itself
nondominating; and consequently all triangles included in $A_c\cup \{c\}$ are dominating. But $A_c$ is not dominating (because there 
is a vertex in $C$ anticomplete to $A_c$), so $|A_c|\le 2$.
\\
\\
(6) {\em $|A_c|=2$.}
\\
\\
Suppose that $|A_c|\le 1$; and hence $|A_c|=1$, since $A_c\cup \{c\}$ is dominating and so is not anticomplete to $S$, since $S\ne \emptyset$.
Let $A_c=\{a\}$; then for the same reason, $a$ is complete to $S$. Also $A_c\cup \{c\}$ is not anticomplete to any vertex in
$A_{c'}$ for $c'\in C\setminus \{c\}$; and so $a$ is complete to $A_{c'}$ for all $c'\in C\setminus \{c\}$. Let $c'\in C\setminus \{c\}$,
and let $a'\in A_{c'}$, and $s\in S$. Since $a',s$ are both adjacent to $a$, and there exists $c''\in C\setminus \{c,c'\}$, and
$G[\{s,a,a',c,c''\}]$ is not a fork, it follows that $a',s$ are adjacent; and so $A_{c'}$ is complete to $S$ for all 
$c'\in C\setminus \{c\}$. If $c_1,c_2\in C\setminus \{c\}$ are distinct, and there exist $a_1\in A_{c_1}$ and $a_2\in A_{c_2}$
nonadjacent, then $G[\{a_1,a_2,a,s,c_1\}]$ is an antifork, a contradiction; so all the sets $A_{c'}\;(c'\in C)$ are complete to each other. Hence they are all homogeneous sets, and so is $S$; so they all have cardinality at most one. But then $G$ is the line graph 
of a bipartite graph, a contradiction. This proves (6).

\bigskip
Let $A_c=\{p,q\}$ say. For each $d\in C\setminus \{c\}$, let $P_{d}$ be the set of vertices in $A_{d}$ adjacent to $p$, and define
$Q_{d}$ similarly. 
\\
\\
(7) {\em For each $d\in C\setminus \{c\}$, $P_d,Q_d$ are disjoint subsets of $A_d$ with union $A_d$.}
\\
\\
Since $A_c\cup \{c\}$ is dominating, it follows that $P_d\cup Q_d=A_d$ for each $d\in C\setminus \{c\}$.
If some $v\in A_d$ is adjacent to both $p,q$, choose $c'\in C\setminus \{c,d\}$; then $G[\{v,p,q,c,c'\}]$ is an antifork, a contradiction.
This proves (7).

\bigskip

Let $S_p,S_q$ be the sets
of vertices in $S$ adjacent to $p,q$ respectively. 
\\
\\
(8) {\em $S_p, S_q$ are disjoint subsets of $S$ with union $S$; and 
$P_d$ is complete to $S_p$, and $Q_d$ is complete to $S_q$, for each $d\in C\setminus \{c\}$.}
\\
\\
Certainly $S_p\cup S_q=S$ since $A_c\cup \{c\}$ is dominating. If $v\in S$
is adjacent to both $p,q$, choose $d\in C\setminus \{c\}$, and then $G[\{v,p,q,c,d\}]$ is an antifork, a contradiction. 
This proves the first assertion.

Let $d\in C\setminus \{c\}$, and suppose
$v\in P_d$ and $s\in S_p$ are nonadjacent. Choose $c'\in C\setminus \{c,d\}$, and then $G[\{v,s,p,c,c'\}]$ is a fork,
a contradiction. Thus $P_d$ is complete to $S_p$, and similarly $Q_d$ is complete to $S_q$, for each $d\in C\setminus \{c\}$.
This proves (8).

\bigskip
Let $P=\{p\}\cup S_p\cup \bigcup_{d\in C\setminus \{c\}}P_d$, and define $Q$ similarly.
\\
\\
(9) {\em $P,Q$ are cliques.}
\\
\\
Suppose $P$ is not a clique, say, and choose $p_1,p_2\in P$, nonadjacent. Thus $p_1,p_2\ne p$. If $p_1\in S$, then by (8), $p_2\in S$,
and $G[\{p_1,p_2,p,c,d\}]$ is a fork where $d\in C\setminus \{c\}$. Thus $p_1,p_2\notin S$, and so we may assume that $p_i\in P_{c_i}$
for $i = 1,2$, where $c_1,c_2,c\in C$ are distinct. But then $G[\{p_1,p_2,p,q,c_1\}]$ is a fork, a contradiction. This proves (9).
\\
\\
(10) {\em For each $d\in C\setminus \{c\}$, $P_d$ is anticomplete to $S_q$, and $Q_d$ is anticomplete to $S_p$.}
\\
\\
Let $v\in P_d$, and suppose $v$ is adjacent to $u\in S_q$. 
Choose $d'\in C\setminus \{c,d\}$; then $G[\{u,v,p,d,d'\}]$
is a fork, a contradiction. This proves (10).
\\
\\
(11) {\em For each $d\in C\setminus \{c\}$, $P_d$ is anticomplete to $Q\setminus Q_d$.}
\\
\\
Suppose that $v\in P_d$ and $u\in Q_{d'}$ are adjacent, where 
$d'\in C\setminus \{c,d\}$. Choose $s\in S$; then $s$ is adjacent to exactly one of $p,q$. If $s$ is adjacent to $p$,
then $s\in P$, and by (9) $s$ is adjacent to $v$; and by (10), $s$ is nonadjacent to $u$. But then $G[\{u,v,s,q,d'\}]$
is a fork, a contradiction. Thus $s$ is adjacent to $q$ and not to $p$, and hence adjacent to $u$ and not to $v$.
But then $G[\{u,v,s,p,d\}]$ is a fork, a contradiction. This proves (11).
\\
\\
(12) {\em Every vertex of $S_p$ has at most one neighbour in $S_q$ and vice versa.}
\\
\\
Suppose $v\in S_p$ is adjacent to $u,w\in S_q$. Then $G[\{v,u,w,q,c\}]$ is an antifork, a contradiction. This proves (12).

\bigskip
From (11), since $P_d$ is complete to $Q_d$, it follows that $P_d$ is a homogeneous set, and so $|P_d|\le 1$, and similarly
$|Q_d|\le 1$.
Let $H$ be the subgraph obtained from $G$ by deleting the edges of the three cliques $C,P,Q$. By (11), every edge of $H$
either has both ends in $A_d\cup \{d\}$ for some $d\in C$, or has one end in $S_p$ and the other in $S_q$;
and hence by (12), every component of $H$
is a clique, and has at most one vertex in common with $C,P$ or $Q$. Consequently $G$ is the line graph of a bipartite graph, 
a contradiction. This proves \ref{nondomtri}.~\bbox

\begin{thm}\label{nomixedtri}
Let $G$ be an uncluttered graph with no nontrivial homogeneous set,
such that $G$ is not the line graph of a bipartite graph. Then either every triangle is dominating, or no clique is dominating.
\end{thm}
\Proof
Let $A$ be union of the vertex sets of all dominating cliques, and let $B$ be the union of the vertex sets of all nondominating
triangles. By \ref{nondomtri}, every vertex of a nondominating triangle only belongs to nondominating cliques, so $A\cap B=\emptyset$.
Suppose there is a nondominating triangle $T$, and a dominating clique $C$. Thus $C\subseteq A$, and $T\subseteq B$.
Every triangle is either dominating or nondominating, and so is a subset of one of $A,B$. Consequently, every vertex of $T$
has at most one neighbour in $C$; and hence exactly one since $C$ is dominating; and since no vertex in $C$ has more than one
neighbour in $T$, it follows that there are three vertices $c_1,c_2,c_3$ such that $t_i$ is adjacent to $c_i$ for $i = 1,2,3$,
where $T=\{t_1,t_2,t_3\}$. By \ref{nondomtri}, $\{c_1,c_2,c_3\}$ is dominating.
Since $T$ is nondominating, there is a vertex $y$ that is anticomplete to $T$; and since $\{c_1,c_2,c_3\}$ is dominating,
we may assume that $y$ is adjacent to $c_1$. If $y$ is nonadjacent to $c_2$ then $G[\{y,c_1,t_1,t_3,c_2\}]$ is a fork, a contradiction.
But if $y$ is adjacent to $c_2$, then since the triangle $\{y,c_1,c_2\}$ has a vertex in $A$, it follows that $y,c_1,c_2\in A$,
and in particular $\{y,c_1,c_2\}$ is dominating; and this is impossible since $t_3$ has no neighbour in  $\{y,c_1,c_2\}$.
This proves \ref{nomixedtri}.~\bbox

The {\em claw} is the complete bipartite graph $K_{1,3}$, and its {\em centre} is its vertex of degree three.
Thus if $T$ is a nondominating triangle in $G$, and $v$ has no neighbour in $T$, then $v$ is a claw centre in $\overline{G}$,
and vice versa. Let us say an {\em anticlaw} is a four-vertex graph whose complement is a claw.
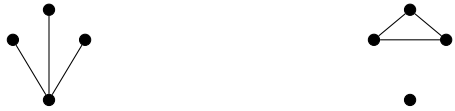
\begin{figure}[H]
\centering

\begin{tikzpicture}[scale=0.8,auto=left]
\tikzstyle{every node}=[inner sep=1.5pt, fill=black,circle,draw]

\node (v1) at (0,0) {};
\node (v2) at (0,1.5) {};
\node (v3) at (-.6,1) {};
\node (v4) at (.6,1) {};

\draw (v1) -- (v2);
\draw (v1) -- (v3);
\draw (v1) -- (v4);

\node (u1) at (6,0) {};
\node (u2) at (6,1.5) {};
\node (u3) at (6-.6,1) {};
\node (u4) at (6+.6, 1) {};

\draw (u2) -- (u3);
\draw (u2) -- (u4);
\draw (u3) -- (u4);

\end{tikzpicture}

\caption{The claw and the anticlaw.} \label{fig:claw}
\end{figure}
We need one more lemma.
\begin{thm}\label{claworanti}
Let $G$ be a graph with no claw or anticlaw. Then for one of $G,\overline{G}$, say $H$, either
\begin{itemize}
\item each component of $H$ is a path or cycle, and hence $H$ is the line graph of a triangle-free graph, or 
\item $|V(H)|\le 9$, and $H$ is the line graph of a bipartite graph.
\end{itemize}
\end{thm}
\Proof
The {\em net} is the graph on six vertices consisting of three pairwise adjacent vertices $t_1,t_2,t_3$, and three more vertices
$s_1,s_2,s_3$, where for $1\le i\le 3$ $s_i$ has degree one and $t_i$ is its unique neighbour. The {\em antinet} is the complement graph
of the net.
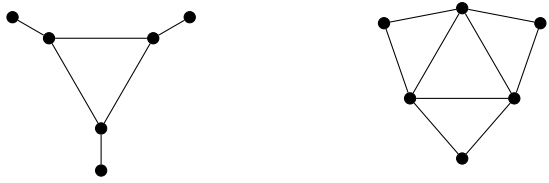
\begin{figure}[H]
\centering

\begin{tikzpicture}[scale=0.8,auto=left]
\tikzstyle{every node}=[inner sep=1.5pt, fill=black,circle,draw]

\def\r{1.5}
\node (v1) at ({sin(0)}, {cos(0)}) {};
\node (v2) at ({sin(120)}, {cos(120)}) {};
\node (v3) at ({sin(240)}, {cos(240)}) {};
\node (v4) at ({\r*sin(60)}, {\r*cos(60)}) {};
\node (v5) at ({\r*sin(180)}, {\r*cos(180)}) {};
\node (v6) at ({\r*sin(300)}, {\r*cos(300)}) {};

\draw (v1) -- (v2);
\draw (v1) -- (v3);
\draw (v2) -- (v3);
\draw (v4) -- (v1);
\draw (v4) -- (v2);
\draw (v5) -- (v2);
\draw (v5) -- (v3);
\draw (v6) -- (v3);
\draw (v6) -- (v1);

\def\r{1.7}
\node (u1) at ({-6+sin(60)}, {cos(60)}) {};
\node (u2) at ({-6+sin(180)}, {cos(180)}) {};
\node (u3) at ({-6+sin(300)}, {cos(300)}) {};
\node (u4) at ({-6+\r*sin(60)}, {\r*cos(60)}) {};
\node (u5) at ({-6+\r*sin(180)}, {\r*cos(180)}) {};
\node (u6) at ({-6+\r*sin(300)}, {\r*cos(300)}) {};

\draw (u1) -- (u2);
\draw (u1) -- (u3);
\draw (u2) -- (u3);
\draw (u4) -- (u1);
\draw (u5) -- (u2);
\draw (u6) -- (u3);

\end{tikzpicture}

\caption{The net and the antinet.} \label{fig:net}
\end{figure}
We begin with:
\\
\\
(1) {\em If $G$ contains a net as an induced subgraph, then $G$ is a net, and so $G$ is the line graph of a bipartite graph.}
\\
\\
Suppose that $s_1,s_2,s_3,t_1,t_2,t_3$ are distinct vertices of $G$ forming a net in the notation above. Let $W=\{s_1,s_2,s_3,t_1,t_2,t_3\}$. 
If $W=V(G)$ the claim holds, so we assume there exists $v\in V(G)\setminus W$.
Since $G[\{v,t_1,t_2,t_3\}]$ is not an anticlaw, $v$ is adjacent to one of $t_1,t_2,t_3$, say $t_1$. 
Since $G[\{v,t_1,s_2,s_3\}]$ is not a claw, $v$ is nonadjacent to one of $s_2,s_3$, say $s_2$. 
Since $G[\{v,s_1,s_2,t_1\}]$ is not an anticlaw, $v$ is nonadjacent to $s_1$. 
Since $G[\{v,s_1,t_1,t_3\}]$ is not a claw, $v$ is adjacent to $t_3$. But then $G[\{v,t_1,s_2,t_3\}]$ is an anticlaw,
a contradiction. This proves (1).

\bigskip
From (1) we may assume $G$ contains no net as an induced subgraph, and by taking complements we may also assume that $G$
contains no antinet. 
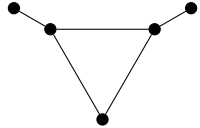
\begin{figure}[H]
\centering

\begin{tikzpicture}[scale=0.8,auto=left]
\tikzstyle{every node}=[inner sep=1.5pt, fill=black,circle,draw]

\def\r{1.7}
\node (u1) at ({6+sin(60)}, {cos(60)}) {};
\node (u2) at ({6+sin(180)}, {cos(180)}) {};
\node (u3) at ({6+sin(300)}, {cos(300)}) {};
\node (u4) at ({6+\r*sin(60)}, {\r*cos(60)}) {};
\node (u6) at ({6+\r*sin(300)}, {\r*cos(300)}) {};

\draw (u1) -- (u2);
\draw (u1) -- (u3);
\draw (u2) -- (u3);
\draw (u4) -- (u1);
\draw (u6) -- (u3);

\end{tikzpicture}

\caption{The bull.} \label{fig:bull}
\end{figure}
The {\em bull} is the graph with five vertices $t_1\ll t_5$, where $t_1\dd t_2\dd t_3\dd t_4$ is an induced path and $t_5$ is adjacent to $t_2,t_3$ and nonadjacent to $t_1,t_4$. Note that the complement of a bull is a bull.
\\
\\
(2) {\em We may assume that $G$ contains a bull as an induced subgraph.}
\\
\\
Suppose first that $G$ has no triangle. Then $G$ has maximum degree at most two, since it has no claw or triangle, and the theorem holds.
So we may assume that $G$ has a triangle, and (by taking complements) $G$ has a stable set of cardinality three.
Choose a triangle $T$, and a set $S$ of three pairwise nonadjacent vertices, with $S\cup T$ minimal.
Let $S=\{s_1,s_2,s_3\}$ and $T=\{t_1,t_2,t_3\}$.
Certainly $|S\cup T|\le 6$; suppose that equality holds. If $t_1$ is adjacent to 
at most one of $s_1,s_2,s_3$, say not to $s_2,s_3$, then $\{t_1,s_2,s_3\}$ is a stable set, contradicting the minimality of $S\cup T$.
So each $t_i$ is adjacent to at least two of $s_1,s_2,s_3$. By the same argument in the complement, each $s_j$ is nonadjacent to at 
least two of $t_1,t_2,t_3$; but this is impossible.
Thus $|S\cup T|\le 5$; and so equality holds, since $|S\cap T|\le 1$. We may assume that $s_3=t_3$. Now each of $s_1,s_2$
is adjacent to one of $t_1,t_2$, since it is not an anticlaw centre; but each of $t_1,t_2$ is adjacent to at most one
of $s_1,s_2$, since it is not a claw centre. Consequently $G[S\cup T]$ is a bull. This proves (2).

\begin{figure}[H]
\centering

\begin{tikzpicture}[scale=0.8,auto=left]
\tikzstyle{every node}=[inner sep=1.5pt, fill=black,circle,draw]
{\scriptsize

\def\r{1.7}
\node (t1) at (-3,0) {};
\node (t2) at (-1, 0) {};
\node (t3) at (1,0) {};
\node (t4) at (3,0) {};
\node (t5) at (0,-1.5) {};

\draw (t5) -- (t2);
\draw (t5) -- (t3);
\draw (t1)--(t2);
\draw (t2) -- (t3);
\draw (t3) -- (t4);

\tikzstyle{every node}=[inner sep=2pt, fill=white,circle,draw]
\node (A1) at (-2,-1.5) {$A_1$};
\node (A2) at (2,-1.5) {$A_2$};
\node (B1) at (-2,1.5) {$B_1$};
\node (B2) at (2,1.5) {$B_2$};
\draw (A1) -- (t1);
\draw (A1) -- (t5);
\draw (A2) -- (t4);
\draw (A2) -- (t5);
\draw (B1) -- (t1);
\draw (B1) -- (t2);
\draw (B1) -- (t4);
\draw (B2) -- (t1);
\draw (B2) -- (t3);
\draw (B2) -- (t4);

\tikzstyle{every node}=[]
\draw (t1) node [left]           {\scriptsize$t_1$};
\draw (t2) node [above right]           {\scriptsize$t_2$};
\draw (t3) node [above left]           {\scriptsize$t_3$};
\draw (t4) node [right]           {\scriptsize$t_4$};
\draw (t5) node [below]           {\scriptsize$t_5$};
}

\end{tikzpicture}

\caption{For step (3) of the proof of \ref{claworanti}.} \label{fig:claworanti}
\end{figure}
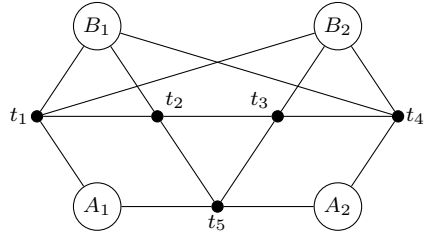

\bigskip
Let $W=\{t_1\ll t_5\}$, and let $G[W]$ be a bull in $G$, with notation as in the definition of a bull. For each vertex
$v\in V(G)\setminus W$, let $W(v)$ denote the set of neighbours of $v$ in $W$. Define $A_1,A_2,B_1,B_2$ by:
\begin{itemize}
\item $A_1$ is the set of $v\in V(G)\setminus W$ with $W(v)=\{t_1,t_5\}$;
\item $A_2$ is the set of $v\in V(G)\setminus W$ with $W(v)=\{t_4,t_5\}$;
\item $B_1$ is the set of $v\in V(G)\setminus W$ with $W(v)=\{t_1,t_2,t_4\}$;
\item $B_2$ is the set of $v\in V(G)\setminus W$ with $W(v)=\{t_1,t_3,t_4\}$.
\end{itemize}
(See figure \ref{fig:claworanti}.)
We claim:
\\
\\
(3) {\em $A_1\cup A_2\cup B_1\cup B_2=V(G)\setminus W$.}
\\
\\
Let $v\in V(G)\setminus W$. Assume first that $t_5\in W(v)$. Since $G[\{v,t_5,t_1,t_4\}]$ is not a claw, one of $t_1,t_4\notin W(v)$,
say $t_1$. Since $G[\{v,t_5,t_3,t_1\}]$ is not an anticlaw, $t_3\notin W(v)$. Since $G[\{v,t_1,t_2,t_3\}]$ is not a claw,
$t_2\notin W(v)$. Since $G[W\cup \{v\}]$ is not a net, $t_4\in W(v)$. But then $v\in A_2$.
Now assume that $t_5\notin W(v)$. We apply the same argument in the complement to deduce that $v\in B_1\cup B_2$. This proves (3).
\\
\\
(4) {\em Each of $A_1,A_2,B_1,B_2$ has cardinality at most one.}
\\
\\
If there exist $u,v\in A_1$, nonadjacent, then $G[\{u,v,t_1,t_2\}]$ is a claw, a contradiction; and if there exist
$u,v\in A_1$ adjacent, then $G[\{u,v,t_1,t_4\}]$ is an anticlaw, a contradiction. Similarly $|A_2|\le 1$; and by taking
complements it follows that $|B_1|,|B_2|\le 1$. This proves (4).
\\
\\
(5) {\em The pairs $(A_1,A_2), (A_1,B_2), (A_2,B_1)$ are complete, and the pairs $(B_1,B_2), (A_1,B_1), (A_2,B_2)$ are anticomplete.}
\\
\\
If there exist $u\in A_1$ and $v\in A_2$, nonadjacent, then $G[\{u,v,t_5,t_2\}]$ is a claw; so $A_1$ is complete to $A_2$.
By taking complement it follows that $B_1$ is anticomplete to $B_2$. If there exists $u\in A_1$ and $v\in B_1$, adjacent, then
$G[\{u,v,t_1,t_3\}]$ is an anticlaw, a contradiction; so $A_1$ is anticomplete to $B_1$, and from the symmetry $A_2$ is anticomplete
to $B_2$. By taking complements it follows that $A_1$ is complete to $B_2$, and $A_2$ is complete to $B_1$. This proves (5).

\bigskip
From (4) and (5) it follows that $|V(G)|\le 9$, and $G$ is an induced subgraph of the line graph of $K_{3,3}$. This proves \ref{claworanti}.~\bbox

We use \ref{claworanti} to prove the next result:
\begin{thm}\label{claw+anticlaw}
Let $G$ be an uncluttered graph with no nontrivial homogeneous set. Then one of $G, \overline{G}$ is the line graph of a triangle-free
graph.
\end{thm}
\Proof Suppose that neither of $G, \overline{G}$ is the line graph of a triangle-free
graph. Suppose first that $G$ contains a claw and an anticlaw. Then since there is a nondominating triangle, it follows from 
\ref{nomixedtri} that no clique is dominating. If $a\in V(G)$ is a claw centre, then it is in at least three 
maximal cliques, and so by \ref{nondomtri}
it is in no triangle (since all triangles are nondominating), and hence the set of neighbours of $a$ is stable. 
By taking complements, if $b\in V(G)$ is an anticlaw centre,
then the set of vertices nonadjacent to $b$ is a clique.
\\
\\
(1) {\em There do not exist a claw centre $a$ and an anticlaw centre $b$ with $a\ne b$.}
\\
\\
Suppose there exist such $a,b$.
By taking complements if necessary, we may assume that $a,b$ are nonadjacent. Since $b$ is an anticlaw centre,
the set of vertices of $G$ nonadjacent to $b$ is a clique $C$
say. Now $|C|\ge 3$ since $b$ is an anticlaw centre; and $a\in C$. Thus $a$ belongs to a triangle, and so belongs to at most two 
maximal cliques by \ref{nondomtri}, contradicting that $a$ is a claw centre. This proves (1).
\\
\\
(2) {\em There do not exist both a claw and an anticlaw in $G$.}
\\
\\
Suppose there is both a claw and an anticlaw. Then by (1), there is a vertex $c$ that is the unique claw centre and the unique 
anticlaw centre. Let $A$ be its set of neighbours and let $B=V(G)\setminus (A\cup \{c\})$. Since $c$ is not in a triangle, $A$ is stable,
and similarly $B$ is a clique. Since $c$ is a claw centre, $|A|\ge 3$, and similarly $|B|\ge 3$. Let $b\in B$. If $b$ has no neighbour
in $A$, then the stable set $\{b,c\}$ is a dominating clique of $\overline{G}$, a contradiction to \ref{nomixedtri} applied to 
$\overline{G}$. If $b$ has at least two nonneighbours in $A$, say $a_1,a_2$, let $b$ be adjacent to $a_3\in A$ and then
$G[\{b,a_1,a_2,a_3,c\}]$ is a fork, a contradiction. But $b$ has at most two neighbours in $A$ since $b$ is not a claw centre;
so $|A|=3$ and $b$ has exactly two neighbours in $A$. By the same argument applied in the complement, $|B|=3$
and every vertex in $A$ has exactly two nonneighbours in $B$, which is impossible by counting edges between $A$ and $B$.
This proves (2).

\bigskip
From (2), and taking complements if necessary, we may assume there is no claw in $G$. 
By \ref{claworanti} we may assume that
there is an 
anticlaw in $G$, that is, there is a nondominating triangle. Consequently every clique is nondominating, by \ref{nomixedtri}.
By \ref{nondomtri} every vertex that is in a triangle is in at most two maximal cliques. But if a vertex $v$ is not in any triangle,
then since there is no claw it follows that $v$ has degree at most two, and so $v$ is in at most two maximal cliques. This proves that
every vertex is in at most two maximal cliques.
Let $C_1\ll C_t$ be the maximal cliques of $G$, and make a graph $H$
with vertex set $\{1\ll t\}$, where distinct $i,j$ are adjacent if $C_i\cap C_j\ne \emptyset$. If $|C_i\cap C_j|\ge 2$ for some distinct $i,j$,
then $C_i\cap C_j$ is a nontrivial homogeneous set of $G$, which is impossible; so $G$ is the line graph of $H$. It remains to show that $H$
is triangle-free. Suppose not; then we may assume that $C_1,C_2,C_3$ pairwise intersect. Choose $v_1\in C_2\cap C_3$, and define $v_2,v_3$ similarly. 
Then $v_{1}\notin C_1$, since every vertex belongs to at most two of $C_1,C_2,C_3$, and similarly $v_{2}\notin C_2$ and 
$v_{3}\notin C_3$. But $v_1,v_2$ are adjacent, because they both belong to $C_3$, and similarly $v_1,v_2,v_3$ are pairwise adjacent;
and hence there is a maximal clique containing all three of $v_1,v_2,v_3$. It is different from $C_1,C_2,C_3$, and so $v_1$
belongs to three different maximal cliques, a contradiction. This proves that $H$ is triangle-free, and so proves \ref{claw+anticlaw}.~\bbox

By combining \ref{claw+anticlaw} and \ref{nohomog}, we deduce our main result, which we restate:
\begin{thm}\label{mainthm2}
Let $G$ be an uncluttered graph. Then either
\begin{itemize}
\item one of $G,\overline{G}$ is disconnected; or
\item one of $G, \overline{G}$ has adjacent simplicial twins; or
\item one of $G,\overline{G}$ is candled; or
\item one of $G,\overline{G}$ is the line graph of a triangle-free graph.
\end{itemize}
\end{thm}
\Proof
If $G$ has a nontrivial homogeneous set, then by \ref{nohomog}
either one of $G,\overline{G}$ is disconnected, or one of $G,\overline{G}$ has
adjacent simplicial twins, or one of $G,\overline{G}$ is candled,
and in each case the 
theorem holds. If $G$ has no nontrivial homogeneous set, then by \ref{claw+anticlaw}, one of $G,\overline{G}$
is the line graph of a triangle-free graph, and again the theorem holds. This proves \ref{mainthm2}.~\bbox

\section{Karthick's question}
We denote the chromatic number of a graph $G$ by $\chi(G)$, and the cardinality of its largest clique by $\omega(G)$.
Let us deduce from \ref{mainthm} a result we stated earlier, that answers a question of Karthick. We restate it:

\begin{thm}\label{chibounded2}
For every uncluttered graph $G$, $\chi(G)\le 2\omega(G)$.
\end{thm}
\Proof
We proceed by induction of $|V(G)|$. We may apply \ref{mainthm}. If $G$ is the disjoint union of two graphs $G_1,G_2$, then
$$\chi(G)=\max(\chi(G_1), \chi(G_2))\le \max(2\omega(G_1), 2\omega(G_2))= 2\omega(G)$$
as required. If $\overline{G}$ is the disjoint union of $\overline{G_1}, \overline{G_2}$, then 
$$\chi(G)=\chi(G_1)+\chi(G_2)\le 2\omega(G_1)+2\omega(G_2)= 2\omega(G)$$
as required.

If $G$ has a simplicial vertex $v$, then we can extend any colouring of $G\setminus \{v\}$ to a 
colouring of $G$ if we have at least $\omega(G)$ colours. Consequently
$$\chi(G)\le \max(\chi(G\setminus \{v\}),\omega(G))\le \max(2\omega(G\setminus \{v\}), \omega(G))\le 2\omega(G)$$
as required. If $G$ has nonadjacent twins $u,v$, then 
$$\chi(G)=\chi(G\setminus \{u\})\le 2\omega(G\setminus \{u\})=2\omega(G)$$
as required. So we may assume that $G$ has no simplicial vertex and no two nonadjacent twins. Consequently
neither $G$ nor $\overline{G}$ has adjacent simplicial twins.

If $G$ is candled, let $Y_1\ll Y_k, Z_1\ll Z_k$ be as in the definition of ``candled''; then any two vertices in $Z_i$ are 
nonadjacent twins, and so each $Z_i$ has cardinality one. But then the vertices in each $Y_i$ are simplicial, a contradiction.
If $\overline{G}$ is candled, again let $Y_1\ll Y_k, Z_1\ll Z_k$ be as in the definition; then any two vertices in $Y_i$ are nonadjacent twins in $G$, so
each $Y_i$ has cardinality one; but then the vertices in each $Z_i$ are simplicial in $G$, a contradiction.

If $G$ is the line graph of a triangle-free graph $H$, then $\chi(G)$ is the edge-chromatic number $\chi'(H)$ of $H$, and
$\omega(G)$ is the maximum degree $\delta(H)$ of $H$. By Vizing's theorem, $\chi'(H)\le \Delta(H)+1$,
so $\chi(G)\le \omega(G)+1\le 2\omega(G)$ (because we can assume that $\omega(G)>0$). Finally, if $\overline{G}$ is the 
line graph of a triangle-free graph $H$, then $\chi(G)$ is the size $\tau(H)$ of the smallest set of vertices of $H$ that meets every edge
of $H$, and $\omega(G)$ is the size $\mu(H)$ of the largest matching in $H$. But $\tau(H)\le 2\mu(H)$, and so
again $\chi(G)\le 2\omega(G)$. This proves \ref{chibounded2}.~\bbox

\section*{Acknowledgements}
We would like to thank T. Karthick for his help with the background to this paper. 
We are also grateful to Jenny Kaufmann for many helpful
discussions, and especially for pointing out to us (in a slightly different
context) the phenomenon of candled graphs. Thanks also to the referee for a very helpful report.

\end{document}